\documentclass[12pt]{amsart}
\usepackage{amsmath,amssymb,amsxtra,amsfonts,latexsym}

\theoremstyle{plain}
\newtheorem{theorem}{Theorem}

\newcommand{\Rcal}{\mathcal{R}}
\newcommand{\C}{\mathbf{C}}
\newcommand{\E}{\mathbf{E}}
\renewcommand{\epsilon}{\varepsilon}
\renewcommand{\phi}{\varphi}
\renewcommand{\ge}{{ }\geqslant{ }}

\newcommand{\Conn}{\mathbf{Conn}}

\title{A Remark on Chern-Weil Theory}
\date{\today}
\author{P. I. Katsylo}
\email{katsylo@gmail.com}

\begin{document}
\begin{abstract}
  We show that Chern-Weil theory for tensor bundles over manifolds is a
  consequence of the existence of natural closed differential forms on total spaces of
  torsion free connections on frame bundles.
\end{abstract}
\maketitle

%%%%%%%%%%%%%%%%%%%%%%%%%%%%%%%%%%%%%
\section{Introduction}
%%%%%%%%%%%%%%%%%%%%%%%%%%%%%%%%%%%%%

Consider $n$-dimensional real manifold $M$ and the frame bundle $\Rcal_M$ over $M$.
We identify bundles and their total spaces. There is the affine bundle
\[
  \pi : \Conn(\Rcal_M) \to M
\]
of torsion free connections on $\Rcal_M$ (sections of bundle $\pi$ are torsion free
connections on $\Rcal_M$).
\begin{theorem}\label{main th}
  There are natural nonzero closed complex valued differential forms
  \[
  \omega_k (M) \in \Omega^{2k} (\Conn(\Rcal_M), \C),
  \qquad k = 1, \ldots , \left[ \frac{n}{2}\right].
  \]
\end{theorem}
In the theorem the adjective \emph{natural} means that \it for every $n$-di\-men\-sio\-nal
real manifold $N$ and open embedding $f : N \to M$ we have
\[
  F^* (\omega_k (M)) = \omega_k (N),
\]
where
\[
  F : \Conn(\Rcal_N) \to \Conn(\Rcal_M)
\]
is the embedding defined by $f$ \cite{curv}. \rm

Suppose that $D$ is a connection on $\Rcal_M$ and $D_0$ is the corresponding
torsion free connection on $\Rcal_M$. So $D_0$ is a differential mapping
\[
  D_0 : M \to \Conn(\Rcal_M).
\]
Consider $2k$-forms
\[
  D_0^*(\omega_k (M)) \in \Omega^{2k}(M, \C), \qquad k = 1, \ldots , \left[ \frac{n}{2}\right].
\]
From theorem \ref{main th} and from the fact that every connection on $\Rcal_M$ can be
deformed to any other one it follows that
\begin{enumerate}
  \item $d (D_0^*(\omega_k (M))) = 0$;
  \item cohomology class $[D_0^*(\omega_k (M))] \in H^{2k}(M, \C)$ does not depend
    on the connection $D$.
\end{enumerate}
In fact, forms $D_0^*(\omega_k (M))$ are exactly $2k$-forms constructed in Chern-Weil
theory \cite{chern-weil}. So theorem \ref{main th} means that Chern-Weil theory for tensor
bundles over manifolds is a consequence of existence of differential forms $\omega_k (M)$.

%%%%%%%%%%%%%%%%%%%%%%%%%%%%%%%%%%%%%%%%
\section{Proof of theorem \ref{main th}}
%%%%%%%%%%%%%%%%%%%%%%%%%%%%%%%%%%%%%%%%

Let $\{ x^i \}$ be local coordinates on an open subset $X$ of the manifold
$M$. The coordinates $\{ x^i \}$ define the local coordinates
$\left\{ x^i, \Gamma_{ij}^k \right\}$ on the open subset $\pi^{-1}(X)$ of the
manifold $\Conn(M)$. On $\pi^{-1}(X)$ consider 1-forms
\[
  \widehat\theta_\beta^\alpha := \Gamma_{i \beta}^\alpha d x^i
\]
and construct the matrix
\[
   \widehat\theta := \left( \widehat\theta_\beta^\alpha \right).
\]
Let $\{ {x'}^i \}$ be local coordinates on an open subset $X$ of the manifold
$M$. By definition of the manifold $\Conn(M)$, on the intersection $X \cap X'$
we have the transition formula
\[
  \widehat\theta' = \left( \frac{\partial x'}{\partial x} \right)^{-1}
  d \left( \frac{\partial x'}{\partial x} \right) +
  \left( \frac{\partial x'}{\partial x} \right)^{-1}
  \widehat\theta \left( \frac{\partial x'}{\partial x} \right).
\]

On $\pi^{-1}(X)$ consider the matrix
\[
  \widehat\Theta := d \widehat\theta + \widehat\theta \wedge \widehat\theta
\]
whose entries are 2-forms on $\pi^{-1}(X)$:
\begin{equation}\label{theta}
  \widehat\Theta_\beta^\alpha = d\Gamma_{i \beta}^\alpha \wedge d x^i \ +
  \ \Gamma_{i k}^\alpha \Gamma_{j \beta}^k dx^i \wedge dx^j.
\end{equation}

We claim that on the intersection $X \cap X'$ we have the transition formula
\begin{equation}\label{tr formula}
  \widehat\Theta' = \left( \frac{\partial x'}{\partial x} \right)^{-1}
  \widehat\Theta \left( \frac{\partial x'}{\partial x} \right).
\end{equation}
We omit an easy check-up of this formula. In fact, it exactly looks like the check-up of
the analogous classical transitional formula.

Define the forms $\omega_k (M)$ by the formula
\[
  1 + \sum_{k \ge 1} \omega_k (M) =
  \det\left(\E + \frac{\sqrt{-1}}{2 \pi} \widehat\Theta \right),
\]
where $\E$ is the identity matrix of size $n \times n$. From \eqref{tr formula} it follows
that forms $\omega_k (M)$ are natural one on the manifold $M$.

Let us prove that forms $\omega_k (M)$ are closed. Let $A \in \Conn(\Rcal_M)$,
$\{ x^i \}$ be local coordinates on a neighbourhood $X \subset M$ of
the point $\pi (A)$, and $\left\{ x^i, \Gamma_{ij}^k \right\}$ be the corresponding
local coordinates on the neighbourhood $\pi^{-1}(X) \subset \Conn(\Rcal_M)$
of the point $A$. Replacing if needed local coordinates $\{ x^i \}$ by new
ones we pass to the case when
\[
  \left. \Gamma_{ij}^k \right|_A = 0.
\]
In this case from \eqref{theta} it follows that
\[
  \left. d \widehat\Theta_\beta^\alpha \right|_A= 0
\]
whence we get
\[
  \left. d \det\left(\E + \frac{\sqrt{-1}}{2 \pi} \Theta \right) \right|_A = 0.
\]
From the above it follows that
\[
  d \det\left(\E + \frac{\sqrt{-1}}{2 \pi} \Theta \right) = 0
\]
and thus forms $\omega(M)$ are closed.

\end{document}